\documentclass[conference,9pt,cmex,intlimits]{sig-alternate}
\usepackage{graphicx}
\usepackage{pstricks}
\usepackage{float}
\usepackage{array}
\usepackage{amsmath}
\usepackage{algorithm,algpseudocode}
\usepackage{gensymb}
\usepackage{framed}
\usepackage[font=small]{caption}
\usepackage{tikz}
\usepackage{pgfplots}

\begin{document}

\title{A Hierarchical Singular Value Decomposition Algorithm for Low Rank Matrices}
\numberofauthors{2}
\author{
  \alignauthor Vinita Vasudevan \\
   \affaddr{Department of Electrical Engineering}\\
       \affaddr{Indian Institute of Technology-Madras}\\
       \affaddr{Chennai-600036, India}\\
        \email{vinita@.iitm.ac.in}
\alignauthor M.Ramakrishna \\
      \affaddr{Department of Aerospace Engineering}\\
       \affaddr{Indian Institute of Technology-Madras}\\
       \affaddr{Chennai-600036, India}\\
      \email{krishna@iitm.ac.in}
}

\maketitle

\begin{abstract}
  Singular value decomposition (SVD) is a widely used technique for dimensionality reduction and computation of basis vectors.
  In many applications, especially in fluid mechanics and image processing the matrices  are dense, but low-rank matrices. In these cases, 
a truncated SVD corresponding to  the most significant singular values is sufficient.
In this paper, we propose a tree based merge-and-truncate algorithm to obtain an approximate truncated SVD of the matrix.
Unlike previous methods, our technique is not limited to ``tall and skinny'' or ``short and fat'' matrices and
it can be used for matrices of arbitrary size.
The matrix is partitioned into blocks and the truncated SVDs of blocks are merged to obtain the final SVD. 
If the matrices are low rank, this algorithm
gives significant speedup over finding the truncated SVD, even when run on a single core. The error
is typically less than 3\%. 
\end{abstract}

\section{Introduction}

Singular Value Decomposition (SVD) is used to obtain basis vectors in a variety of data-driven modelling techniques.
It is a key step in principal component analysis (PCA) (also known as proper orthogonal decomposition (POD)), where
the mean-centered data is arranged as a matrix. 
This is followed by an SVD of the matrix to obtain the basis vectors, which are called the POD modes or eigenfeatures. 
Besides, dimensionality reduction is a key step in many data-driven algorithms such as
facial recognition, latent semantic indexing, collaborative filtering etc., which are used in  the 
evolving data-driven design and modelling algorithms.

In many of  these cases, the matrices are large and very often dense, but  inherently low rank matrices.
Computing the SVD of an $m \times n$ matrix has complexity $O(mn~\textrm{min}(n,m))$.
Since this is super-linear in the size of the data,
it becomes computationally expensive  for large data sets. However, if we have a low rank matrix, we would need only $k$ basis vectors, 
where $k << m,n$. One way of computing the rank $k$ approximation is to compute the SVD of the full matrix
and retain only the $k$ largest singular values and vectors.
It can be shown that this is the best rank $k$ approximation with respect to any unitarily invariant norm. 
The cost of computing this approximation using an SVD followed by truncation  turns out to be expensive, especially if
the matrices are nearer square matrices.
 Moreover, in this ``Bigdata'' era, 
it is entirely possible that the dataset resides on physically different servers and bandwidth and memory constraints on each machine make it impossible
to transfer all the data to a single machine to do the analysis.

Ideally, what is required is a truly distributed algorithm, 
where all the computation is
done in-situ with minimal data transfer and the results of the computation could also reside on several machines.  
This has been attempted in \cite{geist02,bai05,balcan14}. 
These algorithms assume a ``tall and skinny'' data matrix, which is a good assumption for many problems. They partition the matrix row-wise, with
each partition containing a small subset of the rows. If the matrix has $n$
columns, they show that an approximate PCA/SVD can be computed with $O(n^2)$ communication cost. A drawback of the method proposed in \cite{geist02} is that it requires the reconstruction of the low rank approximation of each block of data, followed by accumulation of all these
matrices and an SVD of an $n \times n$ matrix. In \cite{bai05}, they partition the matrix row-wise and perform a hierarchical QR decomposition
by doing a tree-based merge of the $R$ matrices. The SVD of the resultant small $R$ matrix is used to compute the SVD of the full matrix.
In \cite{balcan14}, they once again partition the matrix row-wise and find the SVD of each partition.
This is following by a stacking of the truncated $\Sigma V^T$ of each partition on top of each other and then doing a global SVD.
They also do a randomised SVD of each block to reduce the cost of the block SVDs.
Doing the global SVD can still be quite expensive. Instead of doing a global SVD, a tree-based merging algorithm using truncated  SVDs has been proposed  in \cite{iwen16} to find the SVD of ``short and fat'' matrices.  Here, the partitioning is done column-wise rather than row-wise.

Methods for incremental SVDs have been  proposed in \cite{tufts97,levy00,brand03,brand06,baker12}. 
These algorithms  incrementally compute the SVD of a matrix when new row/columns are added to the matrix. 
They use a combination
of QR and the SVD of a smaller matrix to get the new SVD. The algorithms proposed are essentially sequential algorithms, meant for streaming data.

There are also a number of randomised algorithms, that obtain an approximate truncated SVD in linear time.
A comparison of the performance and accuracy of these algorithms is included in \cite{menon11}. 

In this paper, we propose a hierarchical block based SVD algorithm to obtain a low rank approximation. It combines
the advantages of the methods proposed in 
\cite{geist02,balcan14,tufts97,brand03,iwen16} and is suitable for low rank matrices of arbitrary size.
Unlike \cite{geist02,bai05,iwen16}, our algorithm
is not limited to tall and skinny/short and fat matrices and it is possible to partition the matrix into blocks, both row-wise and column-wise.
The existing algorithms get a runtime improvement when run parallely using several cores. We demonstrate that it is possible to get
speedup even when run on a single core, if partitioned appropriately.
We do tree-based merge of the  truncated SVDs of each block to get the SVD of the matrix
using the incremental SVD method  in \cite{brand03} to find the SVD of the merged blocks.
Each merge is followed by a truncation, where singular values
and vectors that are less than a fraction of the largest singular of the block are discarded. Essentially, our algorithm computes an
approximate truncated SVD of the full matrix using repeated merge-and-truncate (MAT) operations on the blocks. 
 Since we can divide the matrix both
row and column-wise, we do two sets of tree-based MAT operations.
We also propose a iterative method for reducing the error in the approximation.

\section{Incremental and Block-based Algorithms}
\subsection{Algorithms for distributed data-sets }
The distributed SVD algorithms proposed in \cite{bai05,balcan14} are based on the algorithm proposed in \cite{geist02}. It is targeted
to performing a principal component analysis (PCA)  of massive distributed data-sets which reside on several machines and computation is not
possible on a single server. The aim is to minimise communication costs. The assumption here is  that the $m\times n$ matrix $X$, 
is tall and skinny with $m >> n$.
The matrix is partitioned row-wise, with each partition ($X_i$) containing a subset of rows. The steps involved are
\begin{enumerate}
\item Perform an approximate PCA  of $X_i$ locally in each machine. Let $S_i = X_i^T X_i \approx V_i\Sigma_i^2 V_i^T$. This is an $n \times n$ matrix.
\item  The matrices $S_i$  transferred to the central server and summed i.e., $S = \sum_{i=1}^p S_i$. 
\item Perform a PCA of $S$ to get an approximate truncated SVD. 
\end{enumerate}

The algorithm proposed in \cite{bai05}, performs a QR decomposition of each partition instead of an approximate SVD.
The resultant ``R'' matrices from each partition are combined in a tree-based structure to obtain the ``R'' matrix corresponding to the full matrix $X$.
The SVD of the final ``R'' gives the correct $\Sigma V^T$ of $X$.

\subsection{Subspace tracking algorithm}
FAST, proposed by \cite{tufts97}, is an incremental SVD algorithm target-ted to  subspace tracking. Let $X = \begin{bmatrix} X_1 & X_2 & \cdots X_n \end{bmatrix} = U\Sigma V^T$ be an $m\times n$ matrix. If a new column $X_{n+1}$ is added to $X$ and $X_1$ is removed, the goal is to find 
the new SVD incrementally, rather than by re-doing the entire computation.
The idea behind this and similar algorithms is to find the component of $X_{n+1}$ that is orthogonal to the subspace $U$. This is done by subtracting
out the projection of $X_{n+1}$ onto $U$ to obtain the orthogonal component $X_o$ as follows.
\begin{align}
X_o &= X_{n+1} - UU^T X_{n+1} \nonumber \\
q &= \frac{X_o}{||X_o||}
\end{align}
If $X_t =  \begin{bmatrix} X_2 & X_3 & \cdots &  X_{n} \end{bmatrix}$,
\begin{align}
X_{new} &= \begin{bmatrix} X_t & X_{n+1} \end{bmatrix} \nonumber \\
&= \begin{bmatrix} U & q \end{bmatrix} \begin{bmatrix} U^TX_t & U^T X_{n+1} \\ \mathbf{0} & ||X_o|| \end{bmatrix}  \nonumber \\
&=  \begin{bmatrix} U & q \end{bmatrix}  E  
\end{align}
$E$ is an $(n+1) \times n$ matrix and its SVD is inexpensive. If $E = U_E\Sigma_EV_E^T$, the SVD of $X_{new}$ can be written as
\begin{align}
X_{new} &= \begin{bmatrix} U & q \end{bmatrix} U_E\Sigma_EV_E^T \nonumber \\
&= U_N\Sigma_NV_N^T
\end{align}
where $\Sigma_N = \Sigma_E$ and $U_N = \begin{bmatrix} U & q \end{bmatrix} U_E$ and $V_N = V_E$.

\subsection{Online incremental Algorithm}
A generalisation of this algorithm is proposed in \cite{levy00,brand03}. Given $X = U\Sigma V^T$, the author uses a similar technique to find the SVD of
$Y = X + AB^T$. Here $A$ and $B$ can have an arbitrary number of columns and rows respectively. The entries in $A$ and $B$ reflect additions/changes to $X$. Let $Q_AR_A = (I - UU^T)A$ and $Q_BR_B = (I - VV^T)B$ be the QR decomposition of the component of $A$ orthogonal to $U$ and the component of $B$ orthogonal to $V$. This implies
\begin{align}
\begin{bmatrix} U & A \end{bmatrix} &= \begin{bmatrix} U & Q_A\end{bmatrix}\begin{bmatrix} I & U^TA \\ \mathbf{0} & R_A  \end{bmatrix}  \\
\begin{bmatrix} V & B \end{bmatrix} &= \begin{bmatrix} V & Q_B\end{bmatrix}\begin{bmatrix} I & V^TB \\ \mathbf{0} & R_B  \end{bmatrix}
\end{align}
Substituting this in $Y$, we get
\begin{align}
Y &= \begin{bmatrix} U & Q_A\end{bmatrix} \left( \begin{bmatrix} \Sigma & \mathbf{0} \\ \mathbf{0} & \mathbf{0}  \end{bmatrix} + \begin{bmatrix} U^TA \\ R_A \end{bmatrix} \begin{bmatrix} V^TB \\ R_B \end{bmatrix}^T \right) \begin{bmatrix} V^T \\ Q_B^T \end{bmatrix} \nonumber \\
&= \begin{bmatrix} U & Q_A\end{bmatrix} E \begin{bmatrix} V^T \\ Q_B^T \end{bmatrix} \nonumber \\
\end{align}
As in the previous case, the new SVD can be computed using the SVD of the smaller matrix $E$.

Often only a low rank approximation to the subspace is required. $U$ is then an $m \times k$ matrix, where $ k << n,m$. In general, if $p$
additional columns are added, the complexity of the computation is 
$O(2mkp)$ + $O(2mp^2)$ + $O((k+p)^3)$  corresponding to computation of the orthogonal
component, QR decomposition of the orthogonal component and SVD of $E$.

\section{Proposed implementation}
We first start with a simple proof for computing the SVD using either the row or column-wise split, instead of the covariance matrix based
proof in \cite{geist02,iwen16,balcan14}. This proof also extends in a straightforward manner to having simultaneous row and column splits that are needed for
large square matrices. 

Assume the matrix $X$ is an $m \times n$ matrix and is split row-wise into sub-matrices
$X_1$ and $X_2$, with sizes $m_1 \times n$ and $m_2 \times n$ respectively.  
Let $X_1 = U_1\Sigma_1 V_1^T$ and $X_2 = U_2 \Sigma_2 V_2^T$. The SVD of $X = \begin{bmatrix} X_1 \\ X_2 \end{bmatrix}$ can be written as
\begin{align*}
\begin{bmatrix} X_1 \\ X_2 \end{bmatrix} &= \begin{bmatrix} U_1 & \bf{0} \\ \bf{0} & U_2 \end{bmatrix} \begin{bmatrix} \Sigma_1 V_1^T \\ \Sigma_2 V_2^T \end{bmatrix} \nonumber \\
&= \begin{bmatrix} U_1 & \bf{0} \\ \bf{0} & U_2 \end{bmatrix} E
\end{align*}
If $E = \tilde{U} \Sigma V^T$, we get
\begin{align}
\begin{bmatrix} X_1 \\ X_2 \end{bmatrix} &= \begin{bmatrix} U_1 & \bf{0} \\ \bf{0} & U_2 \end{bmatrix} \tilde{U} \Sigma V^T \nonumber \\
 &= U \Sigma V^T
\end{align}
where $U =\begin{bmatrix} U_1 & \bf{0} \\ \bf{0} & U_2 \end{bmatrix} \tilde{U}$ is the product of two orthogonal matrices and hence is also an orthogonal matrix.

Since a low rank approximation to the matrix is required,
we do the merge after truncating the individual SVDs so that only singular values up to $\gamma \sigma_1$ are retained, 
where $\sigma_1$ is the largest singular value of the block. The resultant merged SVD is also truncated using the same criteria.
Hence $X_1 \approx U_{1_k}\Sigma_{1_k} V_{1_k}^T$, $X_2 \approx U_{2_l}\Sigma_{2_l} V_{2_l}^T$
indicating a rank $k$ and $l$ approximation and the result is truncated to a rank-r approximation as follows.
\begin{align}
  \begin{bmatrix} X_1 \\ X_2 \end{bmatrix} & \approx \begin{bmatrix} U_{1_k} & \bf{0} \\ \bf{0} & U_{2_l} \end{bmatrix} \begin{bmatrix} \Sigma_{1_k} V_{1_k}^T \\ \Sigma_{2_l} V_{2_l}^T \end{bmatrix} =  U_r \Sigma_r V_r^T \nonumber  
\end{align}
The algorithm is therefore a merge-and-truncate (MAT) algorithm, rather than a simple merge.

If  $m_1, m_2 < n$, $V_1$ and $V_2$  have dimensions $m_1 \times n$ and $m_2 \times n$. Therefore, the combined  matrix 
$ \begin{bmatrix} \Sigma_1V_1^T \\ \Sigma_2 V_2^T\end{bmatrix}$ is an $m \times n $ matrix. 
By splitting it up and doing three SVDs, we end up actually increasing the computational complexity!
However, if partitioned properly, this process results in a reduction in the computational complexity as will be seen in the next section.

\begin{figure}
\centering
\includegraphics[scale=0.8]{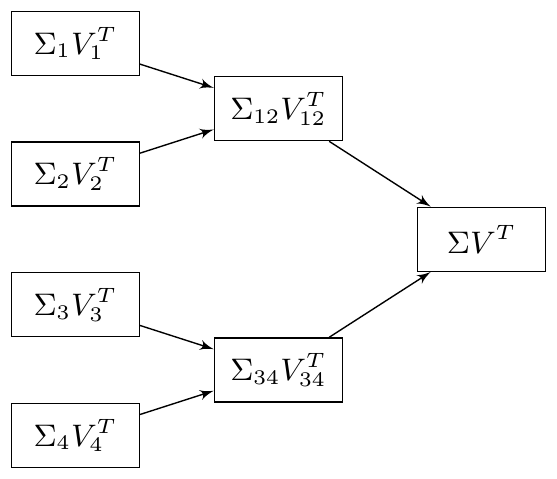}
\caption{Row based partitioning followed by merge of the individual SVDs}
\label{fig:row-merge}
\end{figure}
When there are several partitions, the MAT operations can be done pairwise using a tree 
based algorithm as indicated in Figure \ref{fig:row-merge}.
Note that the ``U'' (left singular-vectors) matrix  is not required for the merge and it need not be propagated.
This is useful since the size of this matrix increases with each merge. 

The matrix can also be split column-wise as $\begin{bmatrix} X_1 & X_2 \end{bmatrix} \approx \begin{bmatrix} U_{1_k}\Sigma_{1_k} V_{_k}^T &  U_{2_l} \Sigma_{2_l} V_{2_l}^T \end{bmatrix}$. In this case, merging can be done as follows.
\begin{align}
\begin{bmatrix} X_1 & X_2 \end{bmatrix} & \approx \begin{bmatrix} U_{1_k} \Sigma_{1_k} & U_{2_l} \Sigma_{2_l} \end{bmatrix} \begin{bmatrix} V_{1_k}^T & \bf{0} \\ \bf{0} & V_{2_l}^T \end{bmatrix} \nonumber \\
 &= U \Sigma \tilde{V}^T  \begin{bmatrix} V_{1_k}^T & \bf{0} \\ \bf{0} & V_{2_l}^T \end{bmatrix} \nonumber \\
 & \approx U_r \Sigma_r V_r^T
\end{align}
Once again, the SVDs of the two blocks are truncated before the merge and the final SVD is also truncated to get a low rank approximation. Also in 
this case, the right singular vectors $V$ are not required for the merge and hence are not propagated through the tree structure.

\begin{figure}
\centering
\includegraphics[scale=0.8]{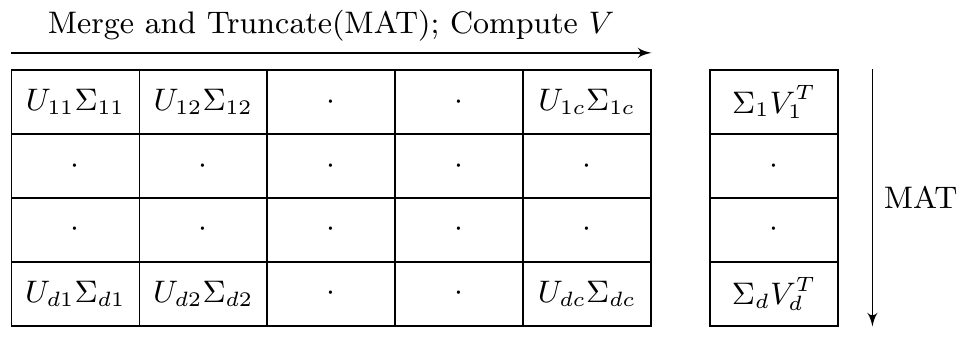}
\caption{Block partitioning consisting of two sets of MAT steps}
\label{fig:bmerge}
\end{figure}

If the matrix is split both row and column-wise, it is easy to see that the blocks can be merged row-wise and column-wise. Let the matrix
$X$ be partitioned as $\begin{bmatrix} X_1 & X_2 \\ X_3 & X_4 \end{bmatrix}$.  A column-wise merge of $X_1, X_2$ and $X_3, X_4$ is followed
by a computation of the  corresponding $V$ matrices. Finally, the two blocks of rows can be merged using $\Sigma V^T$ of each block.
This is illustrated in Figure \ref{fig:bmerge}, for a more general case of a column-wise split into $c$ blocks and a row-wise split of $d$ blocks.
Here, the MAT algorithm is used to first merge column-wise for each set of rows ($X_i$) to obtain $U_i\Sigma_i$. This is followed by computation of
the corresponding $V_i$, which can be obtained by performing an SVD of $U_i^TX_i$.
The MAT algorithm can then be  used to merge  $\Sigma_i V_i$ using a tree structure to give the approximate
singular values and right singular vectors, $\Sigma_r$ and $V_r$.

It is possible to make the merging process more efficient. The MAT algorithm merge requires finding the SVD of matrices of the form
$\begin{bmatrix} \Sigma_{1_k} V_{1_k}^T \\ \Sigma_{2_l} V_{2_l}^T \end{bmatrix}$
and $\begin{bmatrix} U_{1_k} \Sigma_{1_k} & U_{2_l} \Sigma_{2_l} \end{bmatrix}$, which, depending on the block size,  could become expensive.
This can be made more efficient if the orthogonal complement is merged, which can be done using a combination of QR decomposition and an SVD of a smaller matrix using a method similar to \cite{brand03}.
If $X_1 = U_1\Sigma_1V_1^T$ and
$X_2 = U_2 \Sigma_2 V_2^T$ truncated to rank $k$ and $l$ respectively, we can find the component of $U_2$ orthogonal to $U_1$ as $Q = U_2 - U_1U_1^TU_2$. If $Q = U_o R$, 
\begin{align}
\begin{bmatrix} X_1 & X_2 \end{bmatrix} &= \begin{bmatrix} U_1 & U_o \end{bmatrix} \begin{bmatrix} \Sigma_1 & (U_1^T U_2)\Sigma_2 \\ \mathbf{0} & \Sigma_2R \end{bmatrix} \begin{bmatrix} V_1^T & \mathbf{0} \\ \mathbf{0} & V_2^T \end{bmatrix} \nonumber \\
&= \begin{bmatrix} U_1 & U_o \end{bmatrix} E \begin{bmatrix} V_1^T & \mathbf{0} \\ \mathbf{0} & V_2^T \end{bmatrix}
\label{eqn:orthmerge1}
\end{align}
Now $E$ is a much smaller $(k+l) \times (k+l)$ matrix.  If $E = U_E \Sigma_E V_E^T$, we get
\begin{align}
\begin{bmatrix} X_1 & X_2 \end{bmatrix} &= \begin{bmatrix} U_1 & U_o \end{bmatrix} U_E \Sigma_E V_E^T \begin{bmatrix} V_1^T & \mathbf{0} \\ \mathbf{0} & V_2^T \end{bmatrix} \nonumber \\
&= U \Sigma V^T
\label{eqn:orthmerge}
\end{align}
where $U = \begin{bmatrix} U_1 & U_o \end{bmatrix} U_E$, $\Sigma = \Sigma_E$ and $V = V_E \begin{bmatrix} V_1 & \mathbf{0} \\ \mathbf{0} & V_2 \end{bmatrix} $. This is advantageous since we now have to perform an SVD of a much smaller matrix. Although it requires an additional QR decomposition,
the total number of operations  required  is lower for QR than for SVD \cite{chan82,golub}.

Depending on the order of the two sets of MAT operations, we get
the singular values and either the left or right singular vectors.
Assume that after all the MAT operations are over, we have the right singular vectors 
and we need to compute left singular vectors $U_r$. 
Clearly propagation and multiplication of the matrices is not an option. 
Another way to do it is to compute the left singular vectors using the final $\Sigma_r V_r^T$ as $U_r = X V_r \Sigma_r^{-1}$. However, the vectors 
obtained need not be orthogonal, since $\Sigma_r V_r^T$ is approximate due to the truncation at various levels. Instead, we propose the following. 
We know that $V_r$ is an $r \times m$ orthogonal matrix
(assuming a rank $r$ approximation of $X$). Let $X_p = XV_rV_r^T$ be the projection of $X$ onto the $r$ dimensional subspace of the row space. If
$Y = XV_r = U_y \Sigma_y V_y^T$, then $X_p = U_p \Sigma_p V_p^T$ where $U_p = U_y$, $\Sigma_p = \Sigma_y$ and $V_p = V_r V_y$. This is the exact SVD of
the projection of $X$ onto the $r$ dimensional subspace spanned by $V_r$ and an approximation to the best rank-$r$ truncated SVD of $X$.

The steps involved in the algorithm are indicated in Algorithm 1.
The matrix is first partitioned row-wise and each slice is passed to the function DoSVDofColSlices. 
Here, the slice is partitioned column-wise and SVD of each partition is computed. The $U$, $\Sigma$ of each partition are merged using the
function DoMergeOfSlices, which uses the BlockMerge routine. DoSVDofColSlices returns the result of the merger as $U_j$, $\Sigma_j$ for the $j^{\text{th}}$ row slice. This is followed
by computation of the the corresponding $V_j$ and $\Sigma_j$,  as explained.
After the $V_j$ and $\Sigma_j$ of each row $j$ is calculated, they
are merged using the same merge algorithm to obtain $\hat{V}$ and $\hat{\Sigma}$
\begin{algorithm}[htb]
   \caption{Algorithm for block-based SVD}
  \small
  \begin{algorithmic}[1]
    \State $X_{m \times n}$; Input Matrix,  (d, c); Block size
    \State $\gamma$; Truncate parameter used in BlockMerge for MAT. 
    \State 
\Function{DoSVDofBlocks}{$X_{m \times n}$, d, c } 
    \State Nd = round(m/d + 0.5)
    \State $l_X$ = list(), $l_V$ = list(); $l_\Sigma$ = list()
    \State $l_X$ is filled with Nd row slices of X.
    \For {j in range (Nd)}
    \State $U_j$, $\Sigma_j$ = DoSVDOfColSlices($l_X$[j], c)
    \State $V_j$, $\Sigma_j$ = SVD($U_j^T l_X$[j])
    \State $l_V$ += $V_j$; $l_\Sigma$ += $\Sigma_j$
    \EndFor
    \State return $\hat{V},\hat{\Sigma}$ =  DoMergeOfSlices($l_V$, $l_\Sigma$)
    \EndFunction
\State 
    \Function{DoSVDofColSlices}{ X, c }
    \State $l_X$ = list(), $l_U$ = list(), $l_\Sigma$ = list()
    \State Nc = round(n/c + 0.5)
    \State $l_X$ is filled Nc column slices of X.
    \For {j in range(Nc)}
    \State $U_j, \Sigma_j$ = SVD($l_X[j]$)
    \State $l_U$ += $U_j$, $l_\Sigma$ += $\Sigma_j$
    \EndFor
    \State $\hat{U}, \hat{\Sigma}$ = DoMergeOfSlices( $l_U$, $l_\Sigma$)
    \State return $\hat{U}, \hat{\Sigma}$
    \EndFunction       
  \State
  \Function{DoMergeOfSlices}{$l_U$, $l_\Sigma$}
  \State levels = int($\log_2$(len($l_U$)))

    \For {j in range(levels)}
    \State Nl = len($l_U$)
  
  \State $l_{Ut} = l_U$, $l_{\Sigma t} = \Sigma$, $l_U$ = list(), $l_\Sigma$ = list()

  \For {i in range(0, Nl, 2)}
  \State  $U$, $\Sigma$ =  BlockMerge($l_{Ut}[i]$,$l_{\Sigma t}[i]$,$l_{Ut}[i+1]$,$l_{\Sigma t}[i+1]$)
  \State $l_U$ +=  $U_j$, $l_\Sigma$ += $\Sigma_j$
  \EndFor
  \If{ Nl is odd}:
  \State Append last elements of $l_{Ut}$ and $l_{\Sigma t}$ to $l_U$ and $l_\Sigma$
  \EndIf
  \EndFor
  \EndFunction
  \State
  \Function{BlockMerge}{$U_1$, $\Sigma_1$,$U_2$, $\Sigma_2$} 
  \State Use equations (\ref{eqn:orthmerge1}), (\ref{eqn:orthmerge}), and $\gamma$ to do a MAT(QR+SVD)
  \State return $U$, $\Sigma$
  \EndFunction
\end{algorithmic}
\end{algorithm}

Once this is complete, we can follow it up with an iteration to improve accuracy. Algorithm 2 details the steps.
It is essentially equivalent to a power iteration.
In practice, we have observed that only those singular values close to the cut-off value require correction. In all the cases we have seen, two-three iterations proved to be sufficient.

\begin{algorithm}[!t]
\caption{Algorithm for Iterative Improvement}
\small
\begin{algorithmic}[1]
  \State Input: The matrix $X$ and $\hat{V}, \hat{\Sigma}$
  \State Output: $\hat{U}$, $\hat{V}$ and $\hat{\Sigma}$ // Improved accuracy
 
  \State Do:\\
   \hspace{0.1in}   $\tilde{U}_i, \tilde{\Sigma}_i, \tilde{V}_i^T$ = SVD($X \hat{V }$)\\
   \hspace{0.1in}   $\tilde{U}, \tilde{\Sigma}, \tilde{V}^T$  =  SVD($\tilde{U}_i^T X$)\\
   \hspace{0.1in} Error = $||\hat{\Sigma} - \tilde{\Sigma}||_2/||\hat{\Sigma}||_2 $\\
   \hspace{0.1in} $\hat{\Sigma} = \tilde{\Sigma}$\\
   while Error > $\epsilon$:
  \State return $\hat{U} = \tilde{U}_i \tilde{U}, \hat{\Sigma}, \hat{V} = \tilde{V}$ 
\end{algorithmic}
\end{algorithm}

\subsection{Complexity}

In the literature, typically ``tall and skinny'' matrices are partitioned row-wise and ``short and fat'' matrices are partitioned column-wise \cite{balcan14,iwen16}. If this is the case, a runtime speedup can only be obtained when the SVDs of the individual blocks are computed
in parallel fashion with multiple cores. For ``tall and skinny'' matrices, the runtime is $O(mn^2)$. Therefore, if partitioned
column-wise, a runtime improvement can be obtained even on a single machine.

To illustrate this, we do a simplified analysis with all blocks containing the same number of rows and look at the number of floating point operations (flops)
when an  $m \times n$ matrix is partitioned column-wise into $P$ blocks, each containing $s = n/P$ columns. The MAT algorithm
is used to obtain the low rank approximation. Also assume that $m >> n$. the number of flops required for the SVD of the full matrix is approximately $6mn^2 +16n^3$ \cite{bjork,golub}.

The first step is an SVD of each of the $P$ blocks, for which  the number of flops is $P(6ms^2+ 16s^3)$.
Instead of having a truncation based on the magnitude of the singular values, assume that each SVD is truncated to get a rank $k$ matrix.
Therefore, at each level of the binary tree,  the merge cost includes finding the orthogonal complement ($2mk^2$), QR decomposition ($8mk^2$), SVD of a $2k \times 2k$ matrix ($176k^3$), and matrix multiplication to get $U$ is $4mk^2$.
This is done $P-1$ times.
The total number of flops is $P\left[ 6ms^2 + 16 s^3 \right ] + (P-1)\left[ 14mk^2 + 176k^3\right ]$. Since $k \leq s$, the
  total number of flops is  $ < \frac{20mn^2}{P} + \frac{192 n^3}{P^2}$.
  Therefore, we can easily get a speedup.
  
The analysis is similar if the matrix is split into rows and is ``short and fat''. If we have a row and column split, splitting is done 
for the dimension along which the SVD computation is superlinear. Hence for a ``tall and skinny'' matrix, it makes sense to first do the
row-wise split so that the minimum number of rows in a partition is greater than $n$. If this is followed by a column-wise split, it
is possible to  get a speedup.
Note that the speedup can only be obtained if truncation is done before and after the merge.
This in turn is possible only when the original matrix is a low rank matrix.
Also note that this speedup is possible even without parallelization.
Since the first level SVDs and MAT operations at each level of the binary tree can be run independently of each other, there is significant scope for further improvement in the run-time, when run in parallel.

\section{Results}
We have carried out experiments on matrices containing density and velocity data obtained using CFD simulations as well as the FACES dataset \cite{faces}.
All the matrices are dense and low rank. The singular values of the
density dataset as well as the FACES dataset decay more gradually, while the decay is sharp for the velocity data.
For all matrices, we investigated the speedup obtained over performing
a full SVD and then discarding the appropriate number of singular values. The error due to the approximation was measured as $||X_k - \hat{X}_k||_F/||X_k||_F$, where $X_k$ is the rank-k approximation obtained using the full SVD followed by truncation and $\hat{X}_k$ is the rank-k approximation obtained using the various algorithms. 
The code
was written in Python and the SVD routine in LAPACK (available in Scipy), along with the multi-threaded openBLAS library,  was used for algorithm.  We are not presenting results for the parallel version of our code.
The code was run on a core i7 (4 core, 8 threads)  machine running Linux. In all cases, the results reported are the average run times of five runs.
\begin{figure*}
\centering
\includegraphics[width=0.9\linewidth]{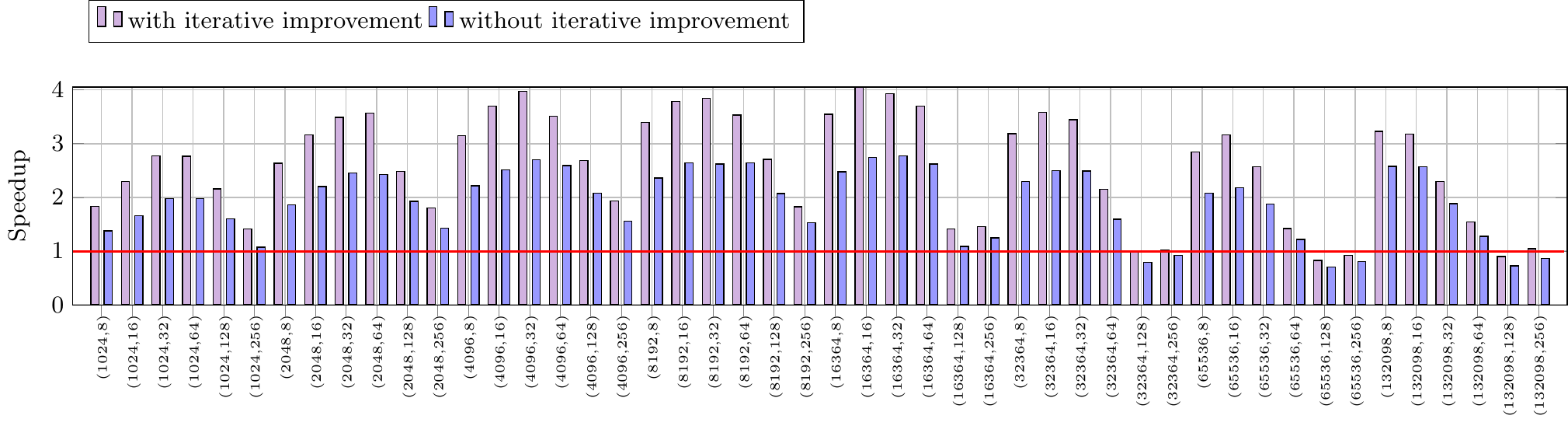}
\caption{Speedup for various block sizes for a velocity dataset obtained from CFD simulations. Merge parameter $\gamma = 10^{-2}$, Iterative improvement convergence criterion $\epsilon = 10^{-3}$. The rank of the approximation obtained was 25-26, except when the number of rows was 132098, when the rank jumped to 32}
\label{fig:suv}
\end{figure*}
\begin{figure*}
\centering
\includegraphics[width=0.9\linewidth]{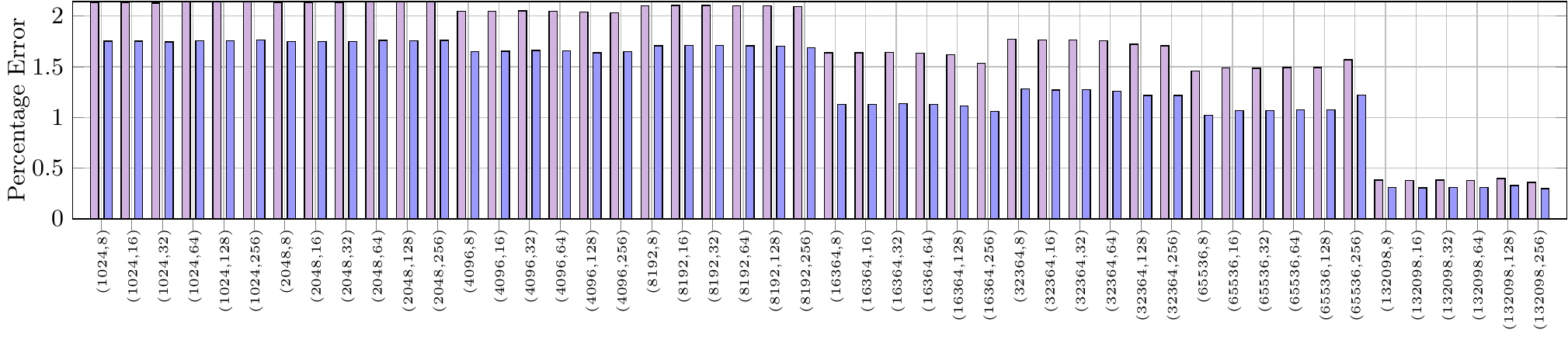}
\caption{Percentage error for various block sizes for a velocity dataset obtained from CFD simulations. Merge parameter $\gamma = 10^{-2}$, Iterative improvement convergence criterion $\epsilon = 10^{-3}$. The rank of the approximation obtained was 25-26}
\label{fig:euv}
\end{figure*}

\begin{table}[H]
  \begin{tabular}{|c|c|c|c|}
    \hline
    & SVD($X$) & SVD($X^T$) & SVD($X^TX$)\\
    \hline
    Run time(s) & 21 & 39.8 & 4.2 \\
    \hline
  \end{tabular}
  \caption{Run time (in seconds) of SVD($X$), SVD($X^T$) and SVD($X^TX$). $X$ is a $132,098 \times 1024$ matrix containing velocity data obtained using CFD simulations. The times reported are the average of 5 runs.}
  \label{tab:rt_svd}
\end{table}
Table \ref{tab:rt_svd} contains the run time in seconds, when an SVD is performed on a $132,098 \times 1024$ matrix containing velocity data.
Note that the SVD computation is distributed across eight threads.
As the table shows, the run times vary significantly depending on whether $X$ or $X^T$ is passed to the LAPACK routine. In general, the run times
were found to be lower when the matrix is ``tall and skinny'' rather than ``short and fat''. Therefore, in all our results for speedup the
run time for SVD$(X)$ was used.
The run time of $X^TX$ includes the time required
for matrix multiplication and time required to compute the left singular vectors.

\begin{table}[H]
  \scriptsize
\begin{tabular}{|c|c|c|}
\hline
 No. of rows & TSQR & Algorithm in \cite{iwen16} \\
\hline
8 & 0.3 & 1.12 (2.14)   \\
\hline
16 & 0.32& 0.94 (2.14)   \\
\hline
128 & 0.37 & 0.47 (2.14)   \\
\hline
512 & 0.43 & 0.2 (2.14)   \\
\hline
2048 & 0.8 &  0.12 (2.03)  \\
  \hline
  4096 & 1 & 0.22 (2) \\
  \hline
  8192 & 1.04 & 0.38 (1.74) \\
  \hline
  16364 & 1.19 & 0.57 (1.57) \\
  \hline
  32768 & 1 & 0.7 (1.53) \\
  \hline
\end{tabular}
\caption{Speedup with respect to the SVD($X$) as a function of the number of rows in each slice using various methods for a $132,098 \times 1024$ matrix. TSQR is the method proposed in \cite{bai05}. The number in brackets is the percentage error.  }
\label{tab:tsqr_mat}
\end{table}

Table \ref{tab:tsqr_mat} contains the speedup with respect to  SVD($X$) of the TSQR algorithm proposed in \cite{bai05} and the algorithm
proposed in \cite{iwen16}. The matrix was split up row-wise and the $R$ matrices in TSQR/ $\Sigma V^T$ matrices in the algorithm proposed in \cite{iwen16}
were merged. In most
cases, we got almost no speedup. Note that, in TSQR, the focus was to reduce communication costs when the matrix is stored in a distributed
fashion in several computers. Also, it gets all the 1024 singular values unlike the algorithm proposed in \cite{iwen16}. 
Partitioning a ``tall and skinny'' matrix row-wise is not expected to  give a speedup unless the SVDs of individual blocks are computed in parallel.
Since the algorithm in \cite{iwen16} gives a low rank approximation, the percentage error is shown in brackets. It is seen to be quite small.

Next, we tried a combined row and column split and used the proposed MAT algorithm for the 132,098 $\times$ 1024 matrix containing velocity data.  
Figures \ref{fig:suv} and \ref{fig:euv} show the speedup and error for the velocity matrix. Without iterative improvement, the algorithm runs about four times faster than doing a full SVD followed by truncation. With iterative improvement, the speedup drops, but it is still above two in most cases. Generally,
it is seen that the speedup is larger when the blocks are ``tall
and skinny''.  The average error is between 1-2\%. As expected, the error drops if iterative improvement is included. It is also apparent that the
error as well as the number of singular values obtained and consequently rank of the approximate matrix is
also relatively independent of block size used. Also significant is the fact that the number of columns in each block can be much lower than the
final number of singular values obtained, so that blocks that have 8 or 16 columns also yield 25 singular vectors.
\begin{figure}[htb]
  \includegraphics{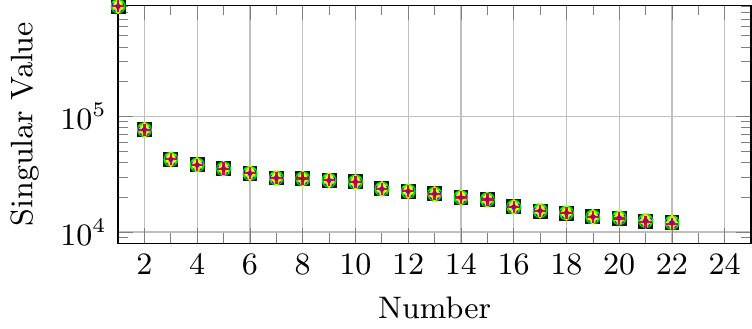}
  \caption{Singular values of the velocity data for various block sizes. The black square represents the  actual singular values}
  \label{fig:sv_uv}
\end{figure}
Figure \ref{fig:sv_uv} shows the singular values obtained for various block sizes for the velocity data as well has the actual singular values of the matrix. It is seen that they match closely, independent of the block size.

\begin{table}[htb]
  \scriptsize
  \begin{tabular}{|c|c|c|c|c|}
    \hline
    $\gamma$ & Max. Speedup & Block size & Error (\%) & Rank\\
    \hline
    0.1 & 9.5 & (132,098$\times$ 8) & 0 & 1\\
    \hline
    0.05 & 5.2 & (132,098 $\times$ 8) & 0.03 & 2\\
    \hline
    0.02 & 3.8 & (132,098 $\times$ 8) & 0.58 & 14 \\
    \hline
    0.01 & 2.8 &(65,536 $\times$ 16)& 1.1 & 26\\
    \hline
    0.005 & 2.33 & (16,364 $\times$ 64) & 0.87 & 53 \\
    \hline
  \end{tabular}
  \caption{Variation with the merge parameter $\gamma$. The table contains the maximum speedup, the block size at which this obtained and the corresponding error and the rank for the velocity matrix.}
  \label{tab:gamma_var}
\end{table}

Table \ref{tab:gamma_var} contains maximum speedup, error and rank of the matrix approximation for various values of the merge parameter
$\gamma$. As expected the speedup decreases as the rank of the approximate matrix increases.
In all cases, the block size for which the maximum speedup was obtained had many more rows than columns.
However, the sensitivity to block size is not very significant, as long as the number of rows is much larger than the number of columns.
For example, with $\gamma = 0.02$, the maximum speedup of 3.8 with a block size of $132 \times 8$.
A block size of $16364 \times 16$, the speedup obtained is 3.4 with an error of 0.65\% and a rank of 11.
It is also seen from the table that larger speedups are obtained when the number of columns in each block is slightly larger than the final
number of singular vectors (rank). As $\gamma$ increases, the rank of the approximated matrix decreases, as expected.
\begin{figure*}
\centering
\includegraphics[width=0.9\linewidth]{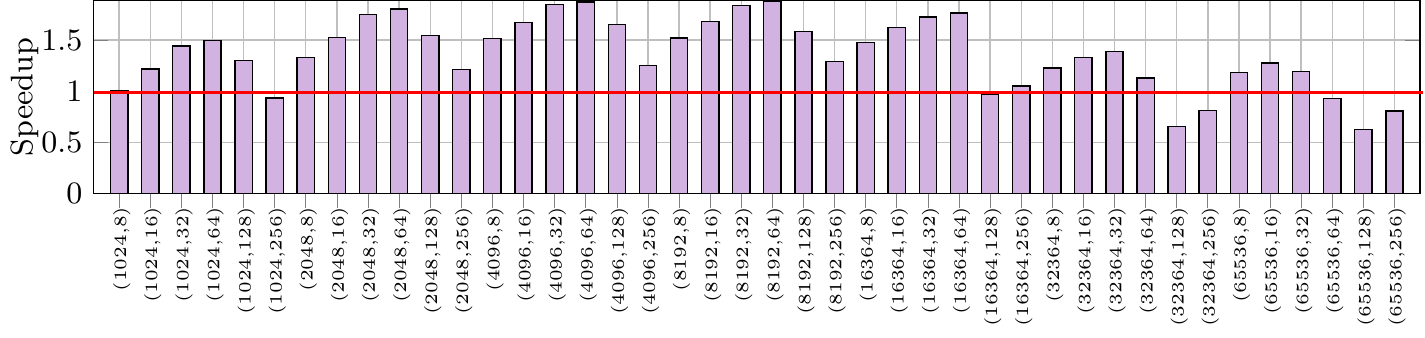}
\caption{Speedup for various block sizes for a density dataset obtained from CFD simulations. Merge parameter $\gamma = 10^{-2}$, Iterative improvement convergence criterion $\epsilon = 10^{-3}$. The rank of the approximation obtained was 65-67}
\label{fig:srho}
\end{figure*}
\begin{figure*}
\centering
\includegraphics[width=0.9\linewidth]{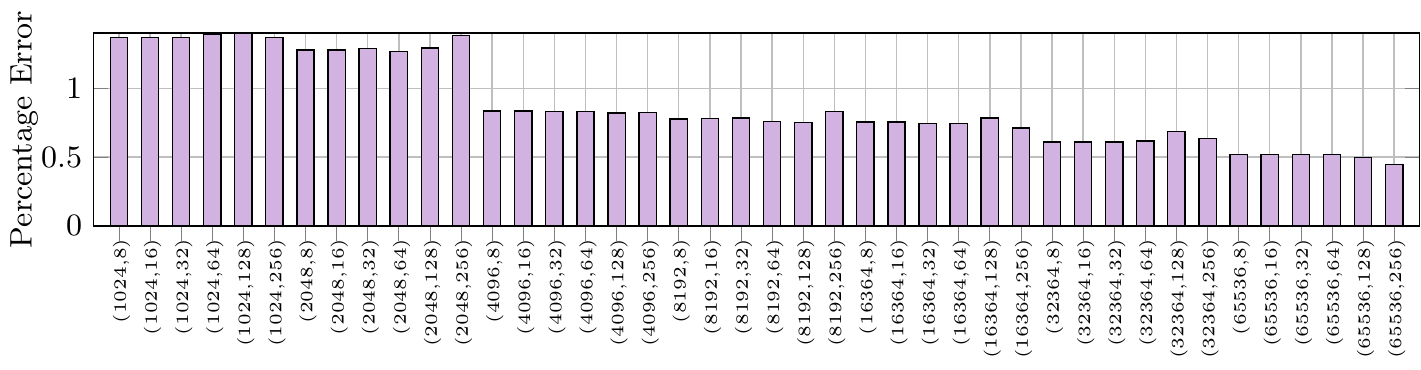}
\caption{Percentage error for various block sizes for a density dataset obtained from CFD simulations. Merge parameter $\gamma = 10^{-2}$, Iterative improvement convergence criterion $\epsilon = 10^{-3}$. The rank of the approximation obtained was 65-67}
\label{fig:erho}
\end{figure*}
\begin{figure*}
\centering
\includegraphics[width=0.9\linewidth]{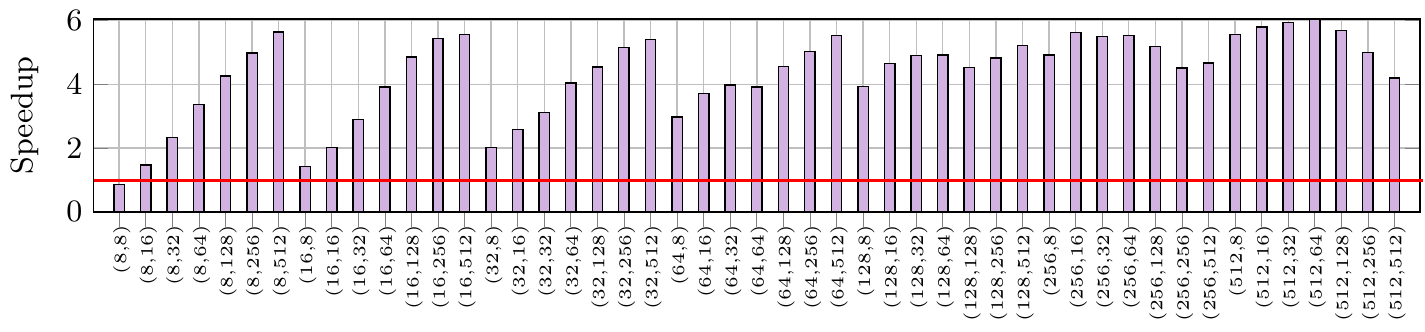}
\caption{Speedup for various block sizes for a velocity dataset with respect to SVD($X^TX$). Merge parameter $\gamma = 10^{-2}$, Iterative improvement convergence criterion $\epsilon = 10^{-3}$. The rank of the approximation obtained was 32-34}
\label{fig:suv_ata}
\end{figure*}
\begin{figure*}
\centering
\includegraphics[width=0.9\linewidth]{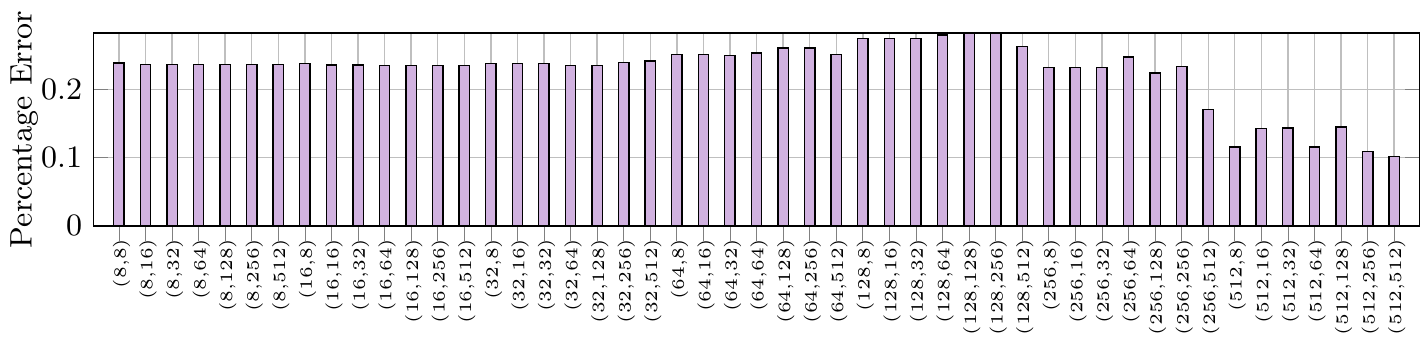}
\caption{Speedup for various block sizes for a velocity dataset with respect to SVD($X^TX$). Merge parameter $\gamma = 10^{-2}$, Iterative improvement convergence criterion $\epsilon = 10^{-3}$. The rank of the approximation obtained was 32-34}
\label{fig:euv_ata}
\end{figure*}

The results for the speedup and error for the density matrix are contained in Figures \ref{fig:srho} and \ref{fig:erho}.
Here, the matrix size is 66049 $\times$ 1024 and the singular values decay more gradually than the velocity matrix.
As mentioned, the speedup is the ratio of the run time of our algorithm and the run time of  doing a full SVD and then truncating.
It is seen that a speedup is obtained in most cases and the largest speedup of about two is obtained when the block size is 8192 $\times$ 64.
The speedup obtained is less than for the velocity data, as the singular values decay more gradually. Therefore,
for the same value of $\gamma$, a larger (65 as opposed to about 25) number of singular values are obtained.

From the results for the density and velocity matrix, it seems like the optimum block size for maximum speedup is when the number of columns is around the number of singular vectors desired and the number of rows in each block is much larger. 

Figures \ref{fig:suv_ata} and \ref{fig:euv_ata} contain the speedup and error with respect to SVD($X^TX$). Here, the matrix $X^TX$ was first
computed and then split into blocks of various sizes. $\gamma$ was set to $10^{-4}$ to allow for comparison with results in Figures \ref{fig:suv} and \ref{fig:euv}. The same value of $\epsilon$ was used. The speedup obtained is significantly larger, as each block is now much smaller.
The error is also significantly lower. We think this is because the rank of the approximated  matrix is larger (between 33-35). As a result,
the large singular values, which typically dominate the error, are approximated better.

\begin{figure}[htb]
  \includegraphics{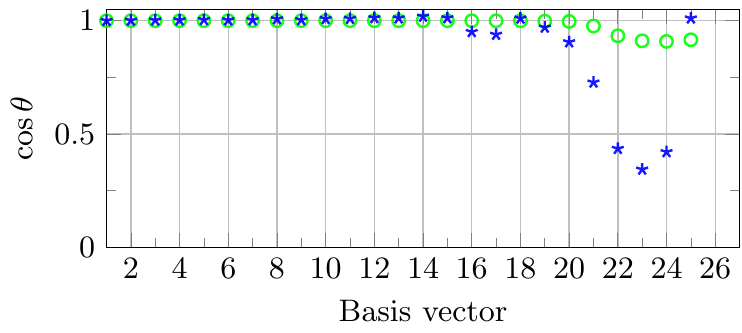}
  \caption{Cosine similarity between actual and approximate left singular vectors with (green circle) and without iterative improvement (blue asterisk). The block size is 16364 $\times$ 16}
  \label{fig:angle_uv}
\end{figure}
\begin{figure}[htb]
  \includegraphics{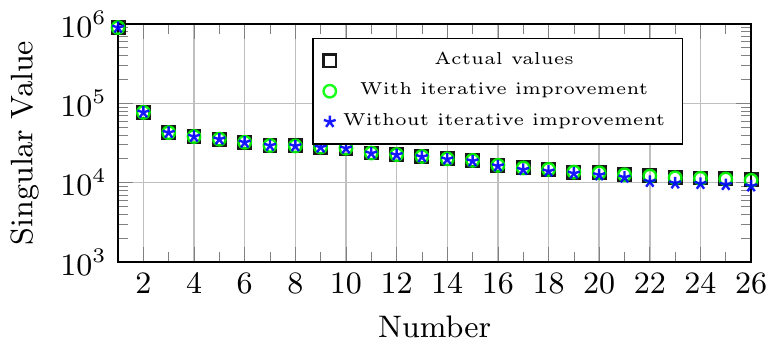}
  \caption{Singular values with and without iterative improvement. The block size is 16364 $\times$ 16}
  \label{fig:sv_uv_new}
\end{figure}
\begin{figure}[!htb]
  \includegraphics{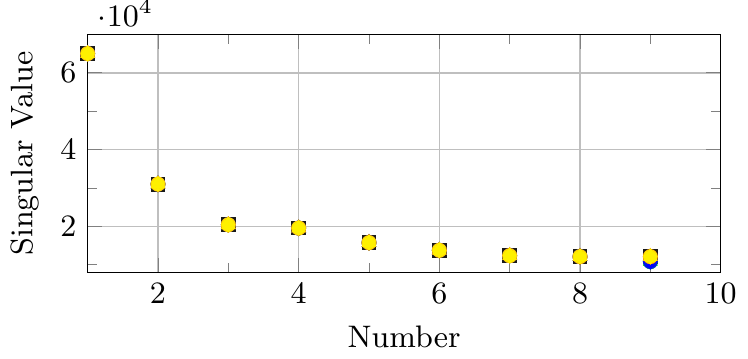}
  \caption{Singular values of the FACES dataset for various block sizes. The black square represents the  actual singular values}
  \label{fig:sv_faces}
\end{figure}
\begin{figure}[htb]
\centering
\includegraphics[width=0.9\linewidth]{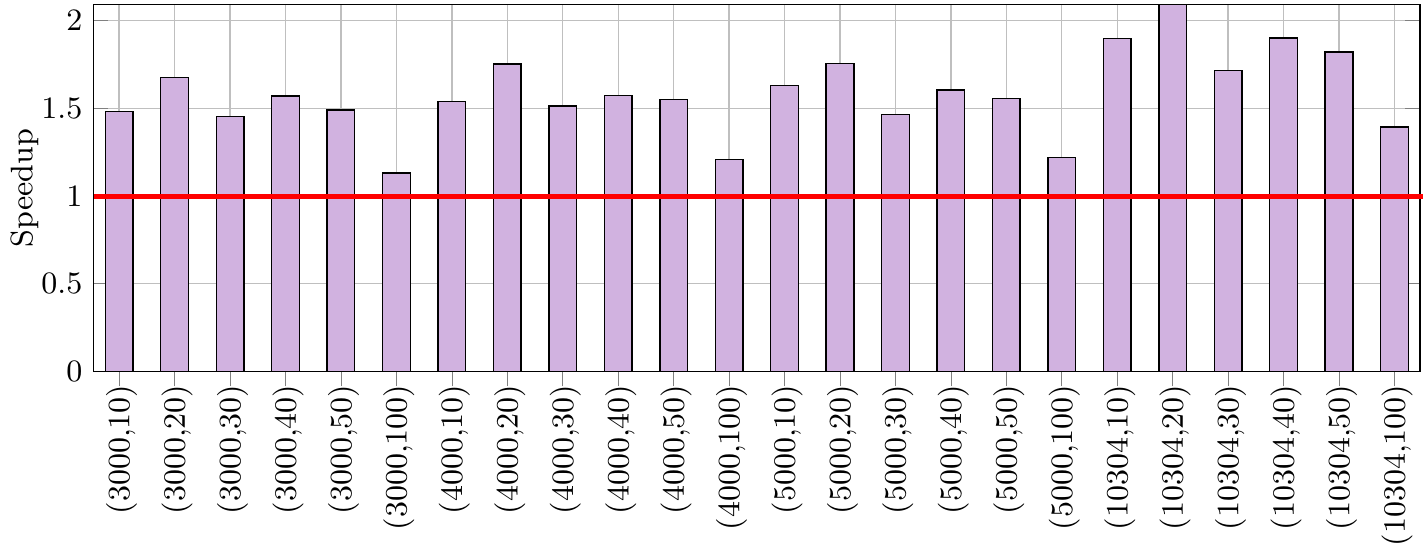}
\caption{Speedup for various block sizes for the FACES dataset. Merge parameter $\gamma = 0.15$, Iterative improvement convergence criterion $\epsilon = 10^{-3}$. The rank of the approximationobtained was 8-9}
\label{fig:sfaces}
\end{figure}
\begin{figure}[ht ]
\centering
\includegraphics[width=0.9\linewidth]{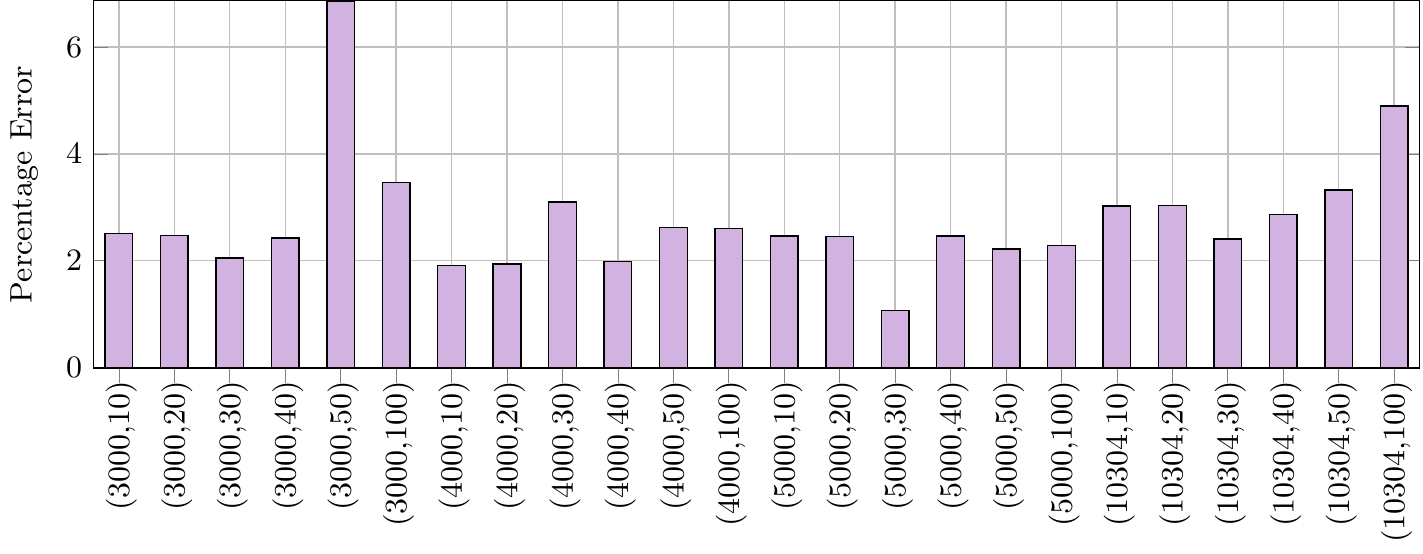}
\caption{Percentage error for various block sizes for the FACES dataset. Merge parameter $\gamma = 0.15$, Iterative improvement convergence criterion $\epsilon = 10^{-3}$. The rank of the approximation obtained was 8-9}
\label{fig:efaces}
\end{figure}

Figures \ref{fig:angle_uv} and \ref{fig:sv_uv_new} show the effect of iterative improvement on the singular values and the angle between
the actual and approximate left-singular vectors, for a block size of 16,364 $\times$ 16. 
It is seen that most of the error is in the last few singular vectors and singular values. The first 19 singular values and vectors match very closely
even without iterative improvement.

We have also done experiments with the ORL database of FACES \cite{faces}. Here, the dimensions of the matrix is 10304 $\times$ 400.
The singular values decay more gradually as shown in Figure \ref{fig:sv_faces}. As seen in the figure, the singular values are approximated well independent of the block size. Figures \ref{fig:sfaces} and \ref{fig:efaces} contain the speedup and percentage error respectively. The
speedup is lower, as it is a smaller matrix to start with.
Typically, most of the error arises due to a poorer
approximation of the last few singular values. Since the decay of singular values is more gradual, the error in these singular values is a larger
percentage of the total error. Therefore, the error in this case is larger than for the CFD datasets.
In fact, the error with a block size of 10304 $\times$ 100 is about 5\% and it
is almost entirely due to the slightly poorer approximation of the last two singular values. 

\section{Conclusions}
In this paper, we have presented a block-based hierarchical merge-and-truncate algorithm to compute a low rank SVD of a matrix. It is suitable for use
in reduced order modelling where the matrices are inherently low rank, but can be quite large and dense. Unlike previous algorithms, it allows for a simultaneous row and column split. We get significant speedup over
doing a full SVD followed by truncation of small singular values, even when run on a single machine. The percentage error is marginal.
The algorithm is very easy to parallelise and can give considerable improvement in the run-time.

Some more work is needed to get error bounds as well as optimal block sizes for various matrices. We have used the SVD routine in LAPACK to obtain the SVD of each block. However, it is possible to use
 any of the randomised algorithms to obtain an approximate SVD of each block. 

\bibliographystyle{ieee}

\end{document}